\documentclass[12pt]{article}

\usepackage{amsmath}
\usepackage{amssymb}
\usepackage{latexcad}

\def\beq{\begin{equation}}
\def\eeq{\end{equation}}
\def\beqq{\begin{eqnarray}}
\def\eeqq{\end{eqnarray}}

\def\Pr{\mathop{{\rm Pr}}\nolimits}
\def\diag{\mathop{{\rm diag}}\nolimits}
\def\E{\mathop{{\rm E}}\nolimits}
\def\dg{\mathop{{\rm dg}}\nolimits}
\def\tr{\mathop{{\rm tr}}\nolimits}

\def\cdc{,\ldots,}

\newcommand{\R}{\mathbb{R}}

\def\withdot#1#2#3{{#1}\kern-#2ex\raise#3ex\hbox{$\centerdot$}\kern#2ex}

\def\si{\sigma}

\def\intercal{\mathop{\scriptscriptstyle\mathrm T}\nolimits}
\def\tra{\intercal}

\def\squareforqed{\hbox{\rlap{$\sqcap$}$\sqcup$}}
\def\qed{\ifmmode\squareforqed\else{\unskip\nobreak\hfil\penalty50\hskip1em\null\nobreak\hfil\squareforqed
\parfillskip=0pt\finalhyphendemerits=0\endgraf}\fi}

\newtheorem{thm}{Theorem}{\bfseries}{\itshape}
{\bfseries}{\itshape}
{\bfseries}{\itshape}
\title{A GRAPH THEORETIC INTERPRETATION OF THE MEAN FIRST PASSAGE TIMES}

\author{Pavel Chebotarev\footnote{E-mail: {\tt pavel4e@gmail.com}}\\
{\normalsize Institute of Control Sciences of the Russian Academy of Sciences}\\
{\normalsize 65 Profsoyuznaya Street, Moscow 117997, Russia}
}

\begin{document}
\maketitle

\unitlength 1.50mm

\begin{abstract}
Let $m_{ij}$ be the mean first passage time from state $i$ to state $j$ in an $n$-state ergodic homogeneous Markov chain
with transition matrix~$T$. Let $G$ be the weighted digraph without loops whose vertex set coincides with the set
of states of the Markov chain and arc weights are equal to the corresponding transition probabilities. We show that
$$
m_{ij}=
\begin{cases}
f_{ij}/q_j,  &\text{if }\;\; i\ne j,\\
1/\tilde q_j,&\text{if }\;\; i=j,
\end{cases}
$$
where $f_{ij}$ is the total weight of 2-tree spanning converging forests in $G$ that have one tree containing~$i$
and the other tree converging to $j$,
$q_j$ is the total weight of spanning trees converging to~$j$ in~$G$, and
$\tilde q_j=q_j/\sum_{k=1}^nq_k$. The result is illustrated by an example.
\bigskip

\bigskip

\end{abstract}

\section{Introduction}

Let $T=(t_{ij})\in\R^{n\times n}$ be the transition matrix of an $n$-state ergodic homogeneous Markov chain
with states $1\cdc n$. Then $T$ is an irreducible stochastic matrix. The {\em mean first passage time from
state $i$ to state $j$\/} is defined as follows:
\beq
\label{defMFPT1}
m_{ij}=\E(F_{ij})=\sum_{k=1}^{\infty}k\Pr(F_{ij}=k),
\eeq
where
\beq
\label{defMFPT2}
F_{ij}=\min\{p\ge1:X_p=j\,|\,X_0=i\}.
\eeq
By \cite[Theorem~3.3]{Meyer75} the matrix $M=(m_{ij})_{n\times n}$ has the following representation:
\beq
\label{Meyer3.3}
M=(I-L^\#+JL_{\dg}^\#)\Pi^{-1},
\eeq
where $L^\#$ is the group inverse of $L$,
\beq
\label{Lapl}
L=I-T,
\eeq
$J$ is the all ones matrix, $L_{\dg}^\#$ is the diagonal matrix obtained by
setting all off-diagonal entries of $L^\#$ to zero, $\Pi=\diag(\pi_1,\ldots,\pi_n)$, and $(\pi_1\cdc\pi_n)=\pi$
is the normalized left Perron vector of $T$, i.e., the row vector in $\R^n$ satisfying
\beq
\label{pi_vect}
\pi T=\pi \:\text{  and  }\: \|\pi\|_1=1.
\eeq

In an entrywise form, (\ref{Meyer3.3}) reads as follows (see, e.g., \cite{CatralNeumannXu05}):
\beq
\label{Meyer3.3entry}
m_{ij}=
\begin{cases}
\pi_j^{-1},                      & \text{ if } i=j,\\
\pi_j^{-1}(L^\#_{jj}-L^\#_{ij}), & \text{ if } i\ne j.
\end{cases}
\eeq

In the following section, we use this formula to derive a graph-theoretic interpretation of the mean first passage times.

\section{A forest expression for the mean first passage times}

Let us say that a weighted digraph $G$ {\em corresponds to the Markov chain with transition matrix~$T$\/} if
the Laplacian matrix $L$
of $G$ satisfies~(\ref{Lapl}). Let the vertex set of $G$ be
$\{1,2\cdc n\}$. Then, by~\eqref{Lapl}, arc $(i,j)$ belongs to the arc set of $G$ whenever [$i\ne j$ and $t_{ij}\ne0$];
the weight of this arc is $t_{ij}$.

Let us recall some graph-theoretic notation. A digraph is {\em weakly connected\/} if the
corresponding undirected graph is connected. A {\em weak component\/} of a digraph $G$ is any maximal weakly
connected subdigraph of~$G$. An {\em in-forest\/} of $G$ is a spanning subgraph of $G$ all of whose weak
components are converging trees (also called {\em in-arborescences}). A {\em converging tree\/} is a weakly
connected digraph in which one vertex, called the {\em root}, has outdegree zero and the remaining vertices
have outdegree one. An in-forest is said to {\em converge to the roots\/} of its converging trees. An in-forest $F$ of
a digraph $G$ is called a {\em maximum in-forest of $G$\/} if $G$ has no in-forest with a greater number of
arcs than in~$F$. The {\em in-forest dimension of a digraph $G$\/} is the number of weak components in any
maximum in-forest. Obviously, every maximum in-forest of $G$ has $n-d$ arcs, where $d$ is the in-forest
dimension of~$G$. A~{\em submaximum in-forest of $G$\/} is an in-forest of $G$ that has $d+1$ weak
components; as a consequence, it has $n-d-1$ arcs.

The weight of a weighted digraph is the product of its arc weights; the weight of any digraph that has no arcs
is~1. The weight of a set of digraphs is the sum of the weights of its members.

By \cite[(iii) of Proposition~15]{CheAga02ap}, for any weighted digraph $G$ and its Laplacian matrix~$L$,
 \beq
 \label{L+viaForests}
 L^\#=\si^{-1}_{n-d}\left(Q_{n-d-1}-\frac{\si_{n-d-1}}{\si_{n-d}}Q_{n-d}\right),
 \eeq
where
$\si_k$ is the total weight of in-forests with $k$ arcs (so that $\si_{n-d}$ and $\si_{n-d-1}$ are the total
weights of maximum and submaximum forests of $G$, respectively), $Q_{k}$ is the matrix whose $(i,j)$ entry
($i,j=1\cdc n$) is the total weight of in-forests that have $k$ arcs and vertex $i$ belonging to the tree
that converges to vertex~$j$.

To obtain a forest representation of the mean first passage times, it suffices to combine
(\ref{Meyer3.3entry}) and~(\ref{L+viaForests}). First, observe that since the Markov chain under
consideration is ergodic, the corresponding digraph $G$ has spanning converging trees. Thus, the in-forest
dimension of $G$ is~1. Consequently, for every $i,j=1\cdc n$, each maximum in-forest converging to $j$ is a
spanning converging tree, which contains~$i$. Therefore, the $(j,j)$ and $(i,j)$ entries of $Q_{n-d}=Q_{n-1}$
are the same.

As a result, substituting \eqref{L+viaForests} in \eqref{Meyer3.3entry} provides the differences of the form
 \beq
 \label{fijdef}
 q^{(n-2)}_{jj}-q^{(n-2)}_{ij}\stackrel{\text{def}}{=}f_{ij}
 \eeq
between the $(j,j)$ and $(i,j)$ entries of the matrix $=(q^{(n-2)}_{ij})=Q_{n-2}=Q_{n-d-1}$. By definition of
$Q_{n-2}$, this difference equals the weight of the set of 2-tree in-forests of $G$ that converge to $j$ and
have $i$ and $j$ in different trees. Thus, substituting (\ref{L+viaForests}) in (\ref{Meyer3.3entry}) yields
 \beq
 \label{eq_mijij}
 m_{ij}=f_{ij}/(\si_{n-1}\,\pi_j) \text{ \;if\; } i\ne j.
 \eeq

Furthermore, we know (for example, from the Markov chain tree theorem
\cite{LeightonRivest83,LeightonRivest86} first obtained in \cite{WentzellFreidlin70a}; see also
\cite{WentzellFreidlin79e,FreidlinWentzell84}) that
 \beq
 \label{eq_MCTT}
 \pi_j=q_j/\si_{n-1},
 \eeq
where $q_j$ is the total weight of trees converging to $j$ in~$G$, so that $\sum_{k=1}^nq_k=\si_{n-1}$. Eqs.
\eqref{eq_mijij}, \eqref{eq_MCTT} and \eqref{Meyer3.3entry} finally provide $m_{ij}=f_{ij}/q_j$ for $i\ne j$
and $m_{jj}=(q_j/\si_{n-1})^{-1}$. We have proved the following theorem.

\begin{thm}
\label{thm_MFPT}
Let $T\in\R^{n\times n}$ be the transition matrix of an $n$-state ergodic homogeneous Markov chain
with states\/ $1\cdc n$. Let $G$ be the weighted digraph without loops whose vertices are\/ $1\cdc n$ and arc weights are
equal to the corresponding transition probabilities in~$T$.
 Then the mean first passage times from state $i$ to state $j$ in the Markov chain can be represented as follows$:$
\beq
\label{MFPT_forest}
m_{ij}=
\begin{cases}
\displaystyle
f_{ij}/q_j,  &\text{if }\;\; i\ne j,\\
\displaystyle
1/\tilde q_j,&\text{if }\;\; i=j,
\end{cases}
\eeq
where $f_{ij}$ is the total weight of\/ $2$-tree in-forests of\/ $G$ that have one tree containing\/ $i$
and the other tree converging to~$j,$
$q_j$ is the total weight of spanning trees converging to~$j$ in $G,$ and
$\tilde q_j=q_j/\sum_{k=1}^nq_k$.
\end{thm}

\noindent
{\bf Remark 1.} If one removes the subsidiary requirement $p\ge1$ from the definition of
mean first passage time (Eq.~(\ref{defMFPT2})), then one has $m_{jj}=0$ and $m_{ij}=f_{ij}/q_j$,
$\;i,j=1\cdc n$, since $f_{jj}=0$ for every~$j$.
\medskip

 \noindent
{\bf Remark 2.} $f_{ij}$ and $q_j$ can be calculated by means of elementary matrix algebra, namely, by the
following recurrence procedure \cite[Proposition~4]{CheAga02ap}. For $k=0,1\cdc$ one has
  \beqq
  \label{req1}
  \si_{k+1}&=&\frac{\tr(LQ_k)}{k+1},\\
  \label{req2}
  Q_{k+1}&=&-LQ_k+\si_{k+1}I,
  \eeqq
where $Q_0=I$.

\section{An example}

In this section, we illustrate Theorem~\ref{thm_MFPT} and Remark~2 by an example.
Let
$$
T=
\left[\begin{array}{rrrr}
  0&   1&    0&    0\\
  0& 0.8&  0.2&    0\\
0.4&   0&  0.2&  0.4\\
  0&   0& 0.25& 0.75\\
\end{array}\right].
$$
By (\ref{Lapl}),
$$
L=
\left[\begin{array}{rrrr}
   1&  -1&     0&    0\\
   0& 0.2&  -0.2&    0\\
-0.4&   0&   0.8& -0.4\\
   0&   0& -0.25& 0.25\\
\end{array}\right].
$$

First, let us obtain the matrix $M$ of the mean first passage times by the direct use of~(\ref{Meyer3.3}).
Finding
\beq
\label{exam_pi0}
\pi=(0.08,0.4,0.2,0.32)
\eeq
and
   \beq
   \label{exam_L+0}
   L^\#=
   \left[\begin{array}{rrrr}
    0.7408&  1.704& -0.448& -1.9968\\
   -0.1792&  2.104& -0.248& -1.6768\\
    0.2208& -0.896&  0.752& -0.0768\\
   -0.0992& -2.496& -0.048&  2.6432\\
   \end{array}\right]
   \eeq
and substituting these in \eqref{Meyer3.3} yields
   \beq
   \label{M0}
   M=
   \left[
   \begin{array}{rrrr}
   12.5 &  1  & 6 & 14.5\\   
   11.5 & 2.5 & 5 & 13.5\\
    6.5 & 7.5 & 5 &  8.5\\
   10.5 &11.5 & 4 &3.125\\
   \end{array}
   \right].
   \eeq
Mention that one method to calculate $L^\#$ is by applying
$$
L^\#=(L+\tilde J)^{-1}-\tilde J,
$$
where
$$
\tilde J=(1,1\cdc\,1)^{\tra}\pi.
$$
(see, e.g., \cite[(i) of Proposition~15]{CheAga02ap}).
\medskip

\begin{figure}[htb]
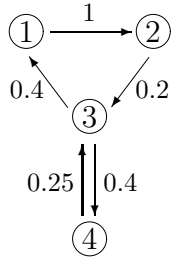
 
\begin{center}
\input gc.lp
\end{center}\vspace{-2em} \caption{The weighted digraph corresponding to the Markov chain.\label{figG}}
\end{figure}

\begin{figure}[htb]
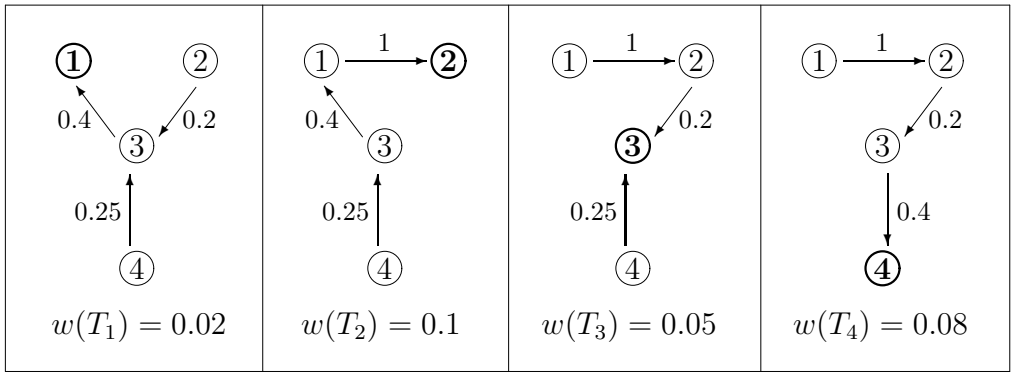
 
\begin{center}
\input t1.lp
\end{center}\vspace{-2em} \caption{The converging trees $T_1,T_2,T_3$, and $T_4$ of $G$.\label{figT}}
\end{figure}

Now let us obtain $M$ by means of Theorem~\ref{thm_MFPT}. The weighted digraph $G$ corresponding to the
Markov chain is shown in Fig.~\ref{figG}. The converging trees of $G$ are shown in Fig.~\ref{figT}, where the roots
are given in a boldface font.
\unitlength 1.50mm

\vspace{-0.1em}
\unitlength 1.50mm

By the definition of $q_i$ given in Theorem~\ref{thm_MFPT} one finds
\beq
\label{q_exam}
q=(q_1,q_2,q_3,q_4)=(0.02,0.1,0.05,0.08).
\eeq
Since $\sum_{k=1}^4q_i=0.25$, by (\ref{q_exam}) it follows that
\beq
\label{tilde_q_exam}
\tilde q=(\tilde q_1,\tilde q_2,\tilde q_3,\tilde q_4)=\frac{q}{\sum_{k=1}^4q_i}=(0.08,0.4,0.2,0.32).
\eeq
This vector coincides with $\pi$, the normalized left Perron vector of~$T$.

\begin{figure}[htb]
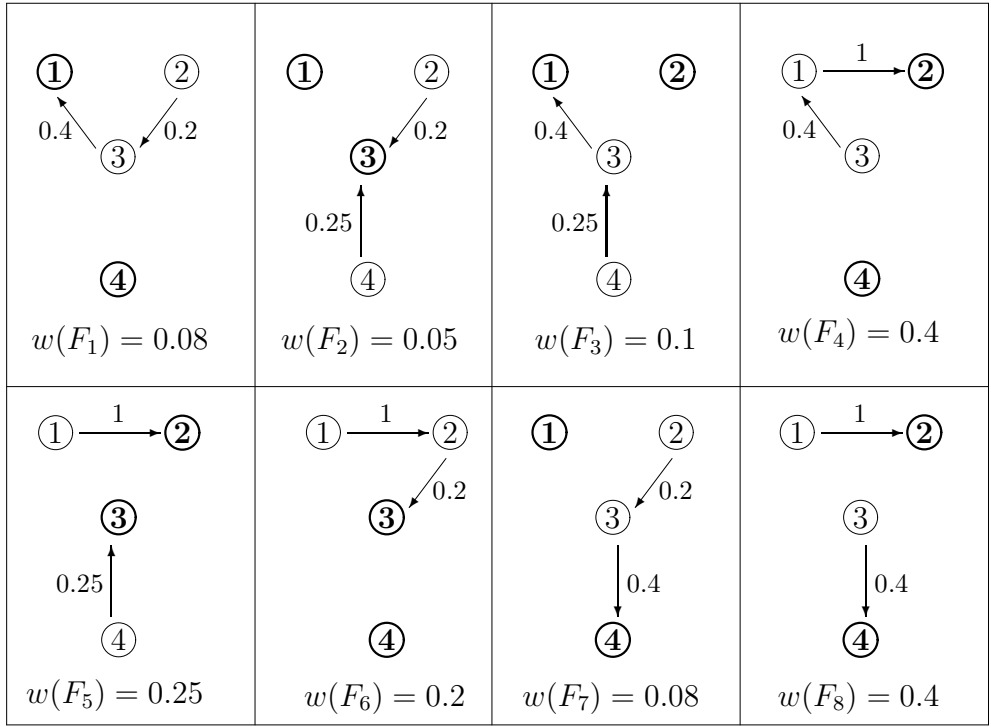
 
\begin{center}
\input f.lp
\end{center}\vspace{-2em} \caption{The 2-tree in-forests $F_1\cdc F_8$ of $G$.\label{figF}}
\end{figure}

The 2-tree in-forests of $G$ are presented in Fig.~\ref{figF}; the roots are shown there in a boldface font.

\unitlength 1.50mm

In Theorem~\ref{thm_MFPT}, $f_{ij}$ is defined as the total weight of 2-tree in-forests of $G$ that have one tree
containing $i$ and the other tree converging to~$j$. Therefore,
   \beqq
   \nonumber
   (f_{ij})
   &=&
   \left[
   \begin{array}{llll}
   0                   &w(\{F_3\})             &w(\{F_2,F_5\}) &w(\{F_1,F_4,F_6,F_7,F_8\})\\
   w(\{F_2,F_3,F_7\})  &0                      &w(\{F_5\})     &w(\{F_1,F_4,F_6,F_8\})    \\
   w(\{F_2,F_7\})      &w(\{F_3,F_5,F_8\})     &0              &w(\{F_1,F_4,F_6\})        \\
   w(\{F_1,F_2,F_7\})  &w(\{F_3,F_4,F_5,F_8\}) &w(\{F_6\})     &0                         \\
   \end{array}
   \right]
   \\
   &=&
   \left[
   \begin{array}{rrrr}
   0    &0.1  &0.3  &1.16\\  
   0.23 &0    &0.25 &1.08\\
   0.13 &0.75 &0    &0.68\\
   0.21 &1.15 &0.2  &0   \\
   \end{array}
   \right],
   \label{exam_fij}
   \eeqq
where $w(A)$ is the weight of a set of digraphs~$A$.

Substituting (\ref{q_exam})--(\ref{exam_fij}) in \eqref{MFPT_forest} yields the matrix of the mean first passage times:
\beq
\label{M1}
M=
\left[
\begin{array}{rrrr}
12.5 &  1  & 6 & 14.5\\   
11.5 & 2.5 & 5 & 13.5\\
 6.5 & 7.5 & 5 &  8.5\\
10.5 &11.5 & 4 &3.125\\
\end{array}
\right].
\eeq

Remark~2 enables one to avoid generating the converging trees and 2-tree in-forests of $G$ shown in
Fig.~\ref{figT} and Fig.~\ref{figF}. Instead, $f_{ij}$ and $q_j$ can be computed by means of the recurrence
procedure (\ref{req1})--(\ref{req2}). Starting with $Q_0=I$, for example under consideration we have:
   \beqq
   \si_1&=&\frac{\tr(LQ_0)}{1}=2.25,\nonumber\\
   Q_1  &=&-LQ_0+\si_1I=
                                 \left[
                                 \begin{array}{rrrr}
                                 1.25&    1&    0&   0\\
                                    0& 2.05&  0.2&   0\\
                                  0.4&    0& 1.45& 0.4\\
                                    0&    0& 0.25&   2\\
                                 \end{array}
                                 \right],\nonumber\\
   \si_2&=&\frac{\tr(LQ_1)}{2}=1.56,\nonumber
   \eeqq
   \beqq
   Q_2  &=&-LQ_1+\si_2I=
                                 \left[
                                 \begin{array}{rrrr}
                                 0.31& 1.05&  0.2&   0\\
                                 0.08& 1.15& 0.25&0.08\\
                                 0.18&  0.4&  0.5&0.48\\
                                  0.1&    0&  0.3&1.16\\
                                 \end{array}
                                 \right],\label{examQ2}\\
   \si_3&=&\frac{\tr(LQ_2)}{3}=0.25,\nonumber\\
   Q_3  &=&-LQ_2+\si_3I=
                                 \left[
                                 \begin{array}{rrrr}
                                 0.02&  0.1& 0.05&0.08\\
                                 0.02&  0.1& 0.05&0.08\\
                                 0.02&  0.1& 0.05&0.08\\
                                 0.02&  0.1& 0.05&0.08\\
                                 \end{array}
                                 \right].\label{examQ3}
   \eeqq

Denoting $Q_2=(q^{(2)}_{ij})_{n\times n}$, by (\ref{fijdef}) we have $f_{ij}=q^{(2)}_{jj}-q^{(2)}_{ij}$,
$i,j=1\cdc n$, hence, by (\ref{examQ2}),
   \beq
   \label{exam_fij1}
   (f_{ij})= \left[
   \begin{array}{rrrr}
   0    &0.1  &0.3  &1.16\\  
   0.23 &0    &0.25 &1.08\\
   0.13 &0.75 &0    &0.68\\
   0.21 &1.15 &0.2  &0\\
   \end{array}
   \right],
   \eeq
which coincides with (\ref{exam_fij}). Eq.~(\ref{examQ3}) yields $q=(0.02,0.1,0.05,0.08)$, which coincides
with (\ref{q_exam}). Thus, we obtain Theorem~\ref{thm_MFPT} provides the matrix~(\ref{M1}) of the mean first
passage times again.



\bibliographystyle{endm}
\bibliography{all2}
\end{document}